\documentclass[11pt]{article}

\nonstopmode

\usepackage[left=1.5cm,right=1.5cm,top=2cm,bottom=2cm]{geometry}
\usepackage[utf8]{inputenc}
\usepackage[french]{babel}
\usepackage[T1]{fontenc}
\usepackage{amsmath}
\usepackage{amsthm}
\usepackage{amsfonts}
\usepackage{amssymb}
\usepackage{color}
\usepackage{multicol}
\usepackage{times}
\usepackage{setspace}
\usepackage{hyperref}
\hypersetup{
colorlinks=true, 
breaklinks=true, 
urlcolor= blue, 
linkcolor= blue, 
citecolor=blue, 
}
\usepackage{graphicx}
\usepackage{lscape}
\usepackage{setspace}

\newcommand{\norm}[1]{\left\lVert #1 \right\rVert}

\newtheorem{thm}{Théorème}[section]
\newtheorem{lem}{Lemme}[section]
\newtheorem{prop}{Proposition}[section]

\newtheorem{con}{Conjecture}[section]

\newtheorem{defn}{Définition}[section]

\newtheorem{ex}{Exemple}[section]

\def\KK{\mathbb K}
\def\QQ{\mathbb Q}

\newcommand{\klammern}[4][]%
{\ifthenelse{\equal{#1}{}}{\left#2}{\csname#1\endcsname#2}%
#4\ifthenelse{\equal{#1}{}}{\right#3}{\csname#1\endcsname#3}}

\newcommand{\abs}[1]{\left\vert#1\right\vert}

\title{On solutions of the Diophantine equation $\mathcal{P}_m-L_n=c$} 

\author{Pagdame TIEBEKABE$^{*,1,2}$, Serge ADONSOU$^3$ $\&$ Ismaïla DIOUF$^1$\\}
\date{}
\begin{document}
\maketitle
\begin{abstract}
\noindent In this article, we determine all the integers $c$ having at least two representations as difference between two linear recurrent sequences. This is a variant of the Pillai's equation. This equation is an exponential Diophantine equation. The proof of our main theorem uses lower bounds for linear forms of logarithms, properties of continued fractions, and a version of the Baker-Davenport reduction method in Diophantine approximation.
\end{abstract}
\begin{itemize}
\item[*]:\small{\textit{ Corresponding author}}
\item[1]:\small{\textit{Université Cheikh Anta Diop (UCAD), Laboratoire d'Algèbre, de Cryptologie, de Géométrie Algébrique et Applications (LACGAA), Dakar, Sénégal}}
\item[2]:\small{\textit{ Université de Kara, Kara-Togo}}
\item[3]:\small{\textit{ African Insitute for Mathematical Sciences (AIMS), Afrique du Sud}}
\end{itemize}

\section{Introduction} 
It is well known that the sequence $\{\mathcal{P}_k\}_{k\geq1}$ of Padovan numbers is defined by
$$
\mathcal{P}_0=\mathcal{P}_1=\mathcal{P}_2=1,\quad \mathcal{P}_{k+3}=\mathcal{P}_{k+1}+\mathcal{P}_k,\quad k\geq 0.
$$
The first Padovan numbers are
$$
1, 1, 1, 2, 2, 3, 4, 5, 7, 9, 12, 16, 21, 28, 37, 49, 65, 86, 114, 151, 200, 265 \ldots
$$
The sequence $\{L_k\}_{k\geq1}$ of Lucas numbers is defined by
$$
L_0=2,\quad L_1=1,\quad L_{k+2}=L_{k+1}+L_k,\quad k\geq0.
$$
The first Lucas numbers are
$$
2, 1, 3, 4, 7, 11, 18, 29, 47, 76, 123, \ldots
$$
In this article, we are interested in the determination of solutions of the Diophantine equation\begin{equation}\label{pleq:main}
\mathcal{P}_m-L_n=c
\end{equation}
for fixed $c$ and $m$, $n$ the unknowns. In particular, we are interested in integers $c$ admitting at least two representations as the difference between a Padovan number and a Lucas number. It is a variant of the equation\begin{equation}\label{pleq2}
a^x-b^y=c,
\end{equation}
in positive integers $(x,y)$ where $a, b, c$ are fixed positive integers. The history of the equation \eqref{pleq2} is very rich and goes back to $1935$.
Subbayya Sivasankaranarayana Pillai ($1901-1950$) is an Indian mathematician specializing in number theory. He has written several articles on perfect powers. A perfect power is a positive integer of the form $a^x$ where $a\geq 1$ and $x\geq 1$ are natural integers.
In $1931$, S.S. Pillai proved in \cite{Pilla} that for all positive integers $a $ and $b$ fixed, both, the number of solutions $(x, y)$ of the Diophantine inequalities $0<a^x - b^y \leq c$ is asymptotically equal to\begin{equation}
\dfrac{(\log c)^2}{2(\log a)(\log b)}.
\end{equation}
when $c$ tends to infinity. It is very interesting to read the $62$ page of this article to see how this result was obtained. This result follows from the attempt to prove that the equation
$$
m^x-n^y=a.
$$
has only a finite number of integral solutions. In this equation, $m$, $n$ and $a$ are fixed. The unknowns are $x$ and $y$. After several years, he repeated this same equation, but this time with $m$, $n$, $x$ and $y$ as unknowns, fixing only $a$.
Research on this equation began with S.S. Pillai in $1931$. In $1936$, A. Herschfeld (\cite{Pl} and \cite{Pil}) continued the research and showed that if $|c|$ is a large enough integer, then the equation
\begin{equation}
2^x-3^y=c
\end{equation}
has at most one solution $(x,y)$ with $x$ and $y$ being positive integers.

This result is no longer true for $|c|$ small enough. By classical methods, Herschfeld demonstrated that only triples of integers $(x,y,c)$ with positive $x$ and $y$ such that $2^x-3^y=c$ is given for $|c|\leq 10$ by:
$$
(2, 1, 1), (1, 1,-1), (3,2,-1), (3,1, 5), (5,3,5), (2,2,5), (4,2,7), (1,2,-7).
$$
So if $x> 5$ or $y> 3$, then $|2^x-3^y|>10$. Proceeding in the same way, he proved that if $x> 8$ or $y> 5$, then $|2^x-3^y|>100$.

S.S. Pillai (\cite{Pl} and \cite{Pil}) extended Herschfeld's results to the more general case of exponential Diophantine equations
\begin{equation}\label{ppeq2}
a^x-b^y=c,
\end{equation}
where $a$ , $b$ and $c$ are nonzero integers fixed with gcd$(a,b)=1$ and $a>b\geq 2$. He showed that there exists a positive integer $c_{0}(a,b)$ such that, for $|c|>c_{0}(a,b)$, this equation has at most one solution. This proof does not give the explicit value of $c_{0}(a,b)$.
In the special case of the Herschfeld equation with $(a,b)=(2,3)$, S.S. Pillai conjectured that $c_{0}(a,b)=13$ and said that the integer $ c$ which has two representations of the form $3^n-2^m$ are the elements of the set $\{-13, -5, 1\}$. This conjecture was solved by R. J. Stroeker and R. Tijdeman in $1982$ by measuring the linear independence of the forms of logarithms of algebraic numbers.
\begin{con} [Conjecture de Pillai]$\\~$
For any integer $k\geq 1$, the Diophantine equation
\begin{equation}
x^n-y^m=k
\end{equation}
admits a finite number of positive integer solutions $(n,m,x,y)$, with $n\geq 2$ and $m\geq 2$.
\end{con}
Since then, several variants of the equation \eqref{ppeq2} have been intensively studied. Recent results related to the equation $H_n-G_n=c$ where $(H_n)_{n\geq 0}$ and $(G_n)_{n\geq 0}$ represent linear recurrent sequences are obtained by M Ddamulira et al in which they solved this type of Pillai equations with Fibonacci numbers and powers of $2$ (see \cite{Ddamulira-Luca-Rakotomalala:2017}), M. Ddamulira et al solved the case with generalized Fibonacci numbers and powers of $2$ (see \cite{Ddamulira-Gomez-Luca:2017}), and Bravo et al solved the case of Tribonacci numbers and powers of $2$ (see \cite{ Bravo-Luca-Yazan:2017}).
We have also solved the case of Padovan numbers and Lucas numbers, by determining the numbers $c$ which have at least two representations as difference of Padovan and Lucas numbers. More simply, we solved the equation $\mathcal{P}_m-L_n=c$ with $m>3$. The articles \cite{16, Paki, Paki1, Paki2, Paki3, Paki4} also discuss variants of the Pillai equation and other Diophantine equations solved by the method of logarithmic linear forms.
The purpose of this article is to prove the following result.
\begin{thm}\label{plth:principal}
The only integers $c$ having at least two representations of the form $\mathcal{P}_m-L_n$ are $$\begin{array}{c}
c \in\{-643, -310, -171, -74, -48, -27, -26, -13, -11, -9, -8, -6, -4, -2, -1, 0, 1, 2, 3,\\ 4, 5, 6, 8, 9, 10, 14, 17, 18, 19, 20, 26, 36, 38, 47, 64, 68, 75, 85, 189, 2864, 58269\}
\end{array}$$
\end{thm}

We organize this article as follows. In the next section, we recall some useful results for the proof of the theorem \ref{plth:principal}. The proof of the theorem \ref{plth:principal} is done in the last section.

\section{Auxiliary results}

\subsection{Some properties of Lucas and Padovan sequences}
\noindent We recall here some properties of Lucas $\{L_k\}$ and Padovan $\{\mathcal{P}_k\}_{k\geq 0}$ sequences which are useful to prove our theorem.\\
The characteristic equation is
$$
x^3-x-1=0,
$$
has roots $\alpha, \beta, \gamma=\overline{\beta}$, where
$$
\alpha=\dfrac{r_1+r_2}{6},\quad \beta=\dfrac{-r_1-r_2+i\sqrt{3}(r_1-r_2)}{12},
$$
and
$$
r_1=\sqrt[3]{108+12\sqrt{69}}\;\text{ and }\ ; r_2=\sqrt[3]{108-12\sqrt{69}}.
$$
Cardan's formulas give for the real root the plastic number or silver number:
$$
\sqrt[3]{\dfrac{1}{2}+\dfrac{\sqrt{69}}{18}}+ \sqrt[3]{\dfrac{1}{2}-\dfrac{\sqrt{69}}{18}}\approx 1,32472.
$$
Also, Binet's formula is
\begin{equation}\label{pleq:Binet_Padovan}
\mathcal{P}_k=a\alpha^k+b\beta^k+c\gamma^k,\ ; \text{ for all } k\geq 0,
\end{equation}
where
\begin{equation}\label{pleq:Formula_a_b_c}
\begin{array}{lll}
a &=& \dfrac{(1-\beta)(1-\gamma)}{(\alpha-\beta)(\alpha-\gamma)}= \dfrac{1+\alpha}{-\alpha^2+3\alpha+1},\vspace{1mm}\\
b &=& \dfrac{(1-\alpha)(1-\gamma)}{(\beta-\alpha)(\beta-\gamma)}=\dfrac{1+\beta}{-\beta^2+3\beta+1},\vspace{1mm}\\
c &=& \dfrac{(1-\alpha)(1-\beta)}{(\gamma-\alpha)(\gamma-\beta)}=\dfrac{1+\gamma}{-\gamma^2+3\gamma+1}=\overline{b}.
\end{array}
\end{equation}
Numerically, we have
\begin{equation}\label{pleq:Numericle_alpha_beta_a_b}
\begin{array}{l}
1.32<\alpha<1.33,\\
0.86<|\beta|=|\gamma|=\alpha^{-1/2}<0.87,\\
0.72 <a<0.73,\\
0.24<|b|=|c|<0.25.
\end{array}
\end{equation}
Using induction, we can show that
\begin{equation}\label{pleq:estimate_Padovan}
\alpha^{k-2}\leq \mathcal{P}_k \leq \alpha^{k-1},
\end{equation}
for all $k\geq 4$.

On the other hand, let $(\delta,\eta)=((1+\sqrt{5})/2,(1-\sqrt{5})/2)$ be the roots of the characteristic equation $ x^2-x-1 = 0$ of the Lucas sequence $\{L_k\}_{k\geq 0}$. Binet's formula for $L_k$

\begin{equation}\label{pleq:Binet-Fibo}
L_k=\delta^k+\eta^k \quad \text{is valid for all } k\geq 0.
\end{equation}
This easily implies that the inequality
\begin{equation}\label{pleq:enca-Fibo}
\delta^{k-1} \leq L_k \leq \delta^{k+1}
\end{equation}
holds for all positive integers $k.$

Now let's discuss the notions of naive height and absolute logarithmic height.

\subsubsection{Algebraic height}
\noindent In this section, we will introduce the notion of algebraic height which is very useful as we will see later. We begin by defining the naive height and then deducing from it the absolute logarithmic height.

\begin{defn} [Naive height]$\\~$
For any algebraic number $\gamma$, we define the height of $\gamma$ by:\begin{align*}
H(\gamma)=\max (|a_d|,\ldots,|a_0|),
\end{align*}
where $f(x)=a_dx^d+\cdots +a_1x+a_0$ is a minimal polynomial of $\gamma$ over $\mathbb{Z}$. $H(\gamma)$ is called the naive height of $\gamma$.\end{defn}

\begin{ex}$\\~$
Let $\alpha$ be an algebraic number:
\begin{itemize}
\item Si $\gamma\in \mathbb{Z}$, $H(\gamma)=|\gamma|$.
\item Si $\gamma\in \mathbb{Q} \left(\text{i.e.} \  \gamma=\dfrac{b}{a}\text{ avec \text{pgcd}}(a,b)=1\right)$, $H(\gamma)=\max\{|a|,|b|\}$,
\end{itemize}
\end{ex}
For any algebraic number $\gamma$, we have the following identity:\begin{equation}
H(\gamma)=|a_d|\prod\limits_{i=1}^d \max \{1,|\gamma_i|\},
\end{equation}
where $\gamma_i$ represent the roots of the minimal polynomial and $f(x)=a_d\prod\limits_{i=1}^d(x-\gamma^{(i)})$ is the minimal polynomial of $\ gamma$. We define in the next subsection, another height deduced from the previous one called absolute logarithmic height. It is the most used.

\begin{defn} [Absolute logarithmic height]$\\~$
For a nonzero algebraic number of degree $d$ over $\mathbb{Q}$ where the minimal polynomial over $\mathbb{Z}$ is \\ $f(x)=a_d\prod\limits_{i=1} ^d(x-\gamma^{(i)})$, we denote by:
\begin{equation}
h(\gamma)=\dfrac{1}{d}\left(\log|a_d|+\sum\limits_{i=1}^d \log\max\{1,|\gamma^i|\}\right)= \dfrac{1}{d} \log  \mathrm M(\gamma),
\end{equation}
the usual absolute logarithmic height of $\gamma$.
\end{defn}
The properties of absolute logarithmic height are as follows:

\begin{prop} [Y. F. Bilu, Y. Bugeaud et M. Mignotte]\label{hauteur}$\smallskip$
\begin{enumerate}
\item Let $\gamma$, $\delta$ be two nonzero algebraic numbers. We have\begin{itemize}
 \item[•]$h(\gamma\delta)\leq h(\gamma)+ h(\delta)$, 
 \item[•]$h(\gamma+\delta)\leq h(\gamma)+h(\delta)+\log 2$.
 \end{itemize}
\item For any algebraic number $\gamma$ and $n\in \mathbb{Z}$ (with $\gamma\neq 0$ and if $n<0$) we have: \\
$h(\gamma^n)=|n|h(\gamma)$.
\end{enumerate}
\end{prop}
More generally, for $\gamma_1,\gamma_2,\cdots,\gamma_n$, $n$ algebraic numbers, we have:
\begin{itemize}
\item[•] $h(\gamma_1\gamma_2\cdots\gamma_n)\leq h(\gamma_1)+h(\gamma_2)+\cdots + h(\gamma_n)$
\item[•] $h(\gamma_1+\gamma_2+\cdots+\gamma_n)\leq h(\gamma_1)+h(\gamma_2)+\cdots+h(\gamma_n)+\log n$.
\end{itemize}

\begin{ex}$\\~$ 
Let $\gamma$ be an algebraic number
\begin{enumerate}
\item If $\gamma$ is the root of $x^2-2x-1$, then $h(\gamma)=\frac{1}{2}(\log\max\{1,|\alpha| \}+\log\max\{1,|\beta|\})$, \\ $\alpha=1+\sqrt{2}$ and $\beta=1-\sqrt{2}$, then $ h(\gamma)=\frac{1}{2}\log\alpha$.
\item Let $h(\gamma)$ be determined with $\gamma=\dfrac{\alpha^n-1}{2\sqrt{2}}$ where $\alpha=1+\sqrt{2}$. From propositions (9.1) and (9.2), we know that $2\sqrt{2}P_n=\alpha^n-\beta^n$ and $Q_n=\alpha^n+\beta^n$.\end{enumerate}
\end{ex}
So,
$$
    \begin{array}{lllll}
        &4\sqrt{2}\gamma+2 &= & 2\sqrt{2}P_n+Q_n\\
        & 4\sqrt{2}\gamma-2\sqrt{2}P_n& =& Q_n-2\\
        &8(4\gamma^2-4P_n\gamma+P_n^2) &= &Q_n^2-4Q_n+4\\
        &8(4\gamma^2-4P_n\gamma)&=&Q_n^2-8P_n^2-4Q_n+4\\
        &8(4\gamma^2-4P_n\gamma)&=&4(-1)^n-4Q_n+4.
    \end{array}
$$
We obtain the minimal polynomial of $\gamma$ divides $8x^2-8P_nx-((-1)^n+1-Q_n)$. Which implies $a(x-\gamma)(x-1+\gamma)$ with $a\in \{1,2,4,8\}$.\begin{align*}
h(\gamma)=\dfrac{1}{a}\left(\log a+\log\max \{1,|\gamma|\}+\log\max\{1,|1-\gamma|\}\right).
\end{align*}

Let us now state the theorems of Stewart, Baker and Wüstholz, before that of Matveev.
\begin{thm} [1993, Baker and Wustholz]$\\~$
If $\Lambda\neq 0$, then
\begin{equation}
|\Lambda|>\exp \left(-(16nd)^{2n+4}.\log A_1...\log A_n.\log B\right)
\end{equation}
with $A_i=\max \{H(\alpha_i),e\}$, for $i=1,\cdots,n$; $B=\max\{|b_1|,\cdots,|b_n|,e\}$ and $d=[\mathbb{Q}(\alpha_1,\cdots,\alpha_n):\mathbb{Q}]$.
\end{thm}

\begin{thm} [A. Baker and G. Wüstholz]$\\~$
Let $K$ be the field of algebraic numbers generated by $\alpha_1,\cdots,\alpha_n$ of degree $d$ over $\mathbb{Q}$. Let $\alpha_1,\cdots,\alpha_n$ $\in K^*$ and $b_1,\cdots,b_n$ $\in\mathbb{Z^*} $.
\end{thm}

Suppose $B^*=\max \{|b_1|,\cdots,|b_n|\}$ and $w=A_1A_2...A_n$ with $|\Lambda|\geq \dfrac{1}{d} \left(\max\{h(\alpha_j),|\log(\alpha_j)|,1\}\right)$ $(1\leq j\leq n)$.

Suppose also that $\Gamma\neq 0$, then;
\begin{center}
$\log(|\Gamma|)>-18(n+1)!n^{n+1}(32d)^{n+2}w\log(2nd)\log(B^*)$.
\end{center}
Let us now state the result of E. Matveev \cite{Matv} which is the most used to solve certain Diophantine equations.
\begin{thm} [E. M. Matveev]\label{theoMatv}$\\~$ 
Let $K$ be an algebraic field of numbers of degree $d$ over $\mathbb{Q}$. If $K\subset\mathbb{R}$, set $\xi=1$ otherwise $\xi=2$. Let $\alpha_1,\ldots,\alpha_n\in \mathbb{K}^*$ and $b_1,\ldots ,b_n\in\mathbb{Z}^*$. Suppose
$$
B^*=\max\{|b_1|,\ldots,|b_n|\}, \  w=A_1A_2\ldots A_n,\  A_j\geq\max\{dh(\alpha_j),|\log(\alpha_j),0.16|\}
$$ 
with ($1\leq j\leq n$) and
$$
\Gamma=b_1\log(\alpha_1)+\cdots +b_n\log(\alpha_n).
$$
If $\Gamma\neq 0$, then $$\log(|\Gamma|)>-C_1(n)d^2w\log(\mathrm{e}d)\log(\mathrm{e}B^* )$$ with $$\ C_1(n)>\min\{\frac{1}{\xi}(0.5\mathrm{e}n)^\xi 30^{n+3}n^{3.5}, 2^{6n+20}\}.$$
\end{thm}

More simply, Y. Bugeaud, M. Mignotte and S. Siksek established the following result.

\begin{thm} [Y. Bugeaud, M. Mignotte, and S. Siksek]$\\~$
Let $n\geq 1$ be an integer. Let $K$ be the field of algebraic numbers of degree $d$. Let $\alpha_1,...,\alpha_n$ be nonzero elements of $K$ and let $b_1,b_2,...,b_n$ be integers,
\begin{center}
$B=\max\{|b_1|,...,|b_n|\}$, 
\end{center}
and 
\begin{center}
$\Lambda=\alpha_1^{b_1}...\alpha_n^{b_n}-1$.
\end{center}
Let $A_1,...,A_n$ be real numbers such that
\begin{center}
$A_j\geq \max\{dh(\alpha_j),|\log (\alpha_j),0.16|\}$, $1\leq j\leq n$.
\end{center}
Assuming $\Lambda\neq 0$, we have:
\begin{center}
$\log|\Lambda|>-3\times 30^{n+4}\times (n+1)^{5.5}\times d^2 \times A_1...A_n(1+\log d)(1+\log nB).$
\end{center}
\end{thm}
If $K$ is real, then
$$
\log|\Lambda|>-1.4\times 30^{n+3}\times (n)^{4.5}\times d^2 \times A_1...A_n(1+\log d)(1+\log B).
$$

Note that for some values of $n$, the lower bound of the logarithm proposed by E.M. Matveev is better (slightly) than that of Baker and Wüstholz.

When $n=2$ and $\alpha_1, \alpha_2$ multiplicatively independent, we have these few results obtained by Laurent, Mignotte, Nesterenko ( \cite{Laurent}, Corollary $2$, pp. 288).

Let in this case $B_1$, $B_2$ be real numbers greater than 1 such that:
$$
\log B_i\geq \max\left\{h(\alpha_i),\dfrac{|\log\alpha_i|}{d},\dfrac{1}{d}\right\}\quad \text{for}\quad i=1,2,
$$
and let's put
$$
b':=\dfrac{|b_1|}{d\log B_2}+\dfrac{|b_2|}{d\log B_1}.
$$
Let's put
$$
\Gamma:=b_1\log \alpha_1+b_2\log \alpha_2.
$$
Note that $\Gamma\neq 0$ because $\alpha_1$ and $\alpha_2$ are multiplicatively independent.

\begin{thm} [Laurent, Mignotte, Nesterenko]$\\~$ \label{theoLaurent}
With the previous notations, let $\alpha_1$, $\alpha_2$ be multiplicatively independent positive numbers, then:
$$
\log |\Gamma|> -24.34 d^4\left( \max\left\{\log b'+0.14, \dfrac{21}{d}, \dfrac{1}{2} \right\}\right)^2\log B_1\log B_2.
$$
\end{thm}
Note that with $\Gamma:=b_1\log \alpha_1+b_2\log\alpha_2$, we have $e^{\Gamma}-1=\Lambda,$ where $\Lambda:=\alpha_1^{b_1} \cdots\alpha_n^{b_n}$ in case $n=2.$

\section{Reduction method}
\noindent During calculations, we get upper bounds on our variables which are too large, so we have to reduce them. To do this, we use some results from the theory of continued fractions. Concerning the treatment of homogeneous linear forms in two integer variables, we use the well-known method of the classical result in the theory of Diophantine approximation.

\begin{lem}[Legendre]\label{lemLegendre}$\\~$
Let $\tau$ be an irrational number, $\dfrac{p_0}{q_0}, \dfrac{p_1}{q_1},\dfrac{p_2}{q_2},\ldots$ all the convergents of the continued fraction of $\ tau$, and $M$ a positive integer. Let $N$ be a positive integer such that $q_N>M.$ Then setting $a(M):=\{a_i:i=0,1,2,\ldots,N\},$ the inequality
$$
\left|\tau-\dfrac{r}{s}\right|> \dfrac{1}{(a(M)+2)s^2},
$$
is valid for all pairs $(r, s)$ of positive integers with $0<s<M.$
\end{lem}

For a non-homogeneous linear form with two integer variables, we use a slight variation of a result due to Dujella and Peth\H o (\cite{Dujella}, Lemma 5a). The proof is almost identical to that of the corresponding result in \cite{Dujella}.
For a real number $X$, we write $\norm{X} :=\min\{\mid X-n\mid :n\in \mathbb{Z}\}$ for the distance from $X$ to the nearest integer.

\begin{lem} [Dujella, Peth\H o]\label{lemDujella}$\\~$
Let $M$ be a positive integer, $\dfrac{p}{q}$ a convergent of the continued fraction of the irrational number $\tau$ such that $q> 6M$, and $A, B, \mu$ be numbers algebraic such as $A>0$ and $B>1.$ Also, $\varepsilon:=\norm{\mu q}-M\norm{\tau q}$. If $\varepsilon>0$, then the following inequality:
$$
0<|u\tau-v+\mu|<AB^{-w},
$$
does not admit an integer solution $u$, $v$ and $w$ with
$$
 u\leq M\quad \text{et}\quad w\geq \dfrac{\log(Aq/\varepsilon)}{\log B}.
$$
\end{lem}
On various occasions we need to find a lower bound for linear forms of logarithms with bounded integer coefficients in three and four variables. In this case, we use the Lenstra-Lenstra-Lovász basic lattice reduction algorithm (LLL-algorithm) which we describe below. Let $\tau_1,\tau_2,\ldots,\tau_t \in \mathbb{R}$ and the linear form
\begin{equation}
x_1\tau_1+x_2\tau_2+\cdots+x_t\tau_t\quad \text{with}\quad |x_i|\leq X_i.
\end{equation}

We set $X:=\max \{X_i\}$, $C> (tX)^t$ and consider the entire lattice $\Omega$ generated by:
$$
b_j:=e_j + \lfloor c \tau_j \rceil \quad \text{for} \quad 1\leq j\leq t-1 \quad \text{and} \quad b_t:=\lfloor C\tau_t\rceil e_t,
$$
where $C$ is a sufficiently large positive constant.

\begin{lem} [LLL-algorithme]$\\~$
Let $X_1, X_2, \cdots, X_t$ be positive integers such that $X:=\max \{X_i\}$ and $C> (tX)^t$ is a sufficiently large fixed positive constant. With the above notations on the lattice $\Omega$, we consider a reduced basis ${b_i}$ to $\Omega$ and its associated Gram-Schmidt orthogonalization basis $\{b_i^*\}$. We fix
$$\displaystyle
c_1:=\max_{1\leq i\leq t} \dfrac{\norm{b_1}}{\norm{b_i^*}}, \quad \theta:=\dfrac{\norm{b_1}}{c_1}, \quad Q:=\sum _{i=1}^{t-1} X_i^2, \quad \text{et}\quad R:=\dfrac{1}{2}\left(1+ \sum _{i=1}^{t} X_i \right). 
$$
If the integers $x_i$ are such that $|x_i|\leq X_i$, for $1\leq i\leq t$ and $\theta^2\geq Q+R^2$, then we have
$$
\left | \sum_{i=1}^{t}x_i\tau_i \right|\geq \dfrac{\sqrt{\theta^2-Q}-R}{C}.
$$
\end{lem}
For proof and further details, we refer the reader to Cohen's book. (Proposition $2.3.20$ in \cite{Cohen}, pp. $58-63$).
\section{Main result}

\noindent Suppose there are positive integers $n,m,n_1,m_1$ such that $(n,m)\neq (n_1,m_1),$ and
$$
L_n-\mathcal{P}_m=L_{n_1}-\mathcal{P}_{m_1}.
$$
Due to symmetry, we can assume that $m \geq m_1$. If $m=m_1$, then $L_n=L_{n_1}$, thus $(n, m)=(n_1,m_1)$, contradicting our hypothesis. Thus, $m > m_1$. Seen that
\begin{equation}\label{pleq:Pell_2}
L_n-L_{n_1}=P_m-P_{m_1},
\end{equation}
and the right member is positive, we obtain that the left member is also positive and therefore $n>n_1$. Thus, $n\geq 2$ and $n_1\geq 1$. Using Binet's formulas \eqref{pleq:Binet-Fibo} and \eqref{pleq:Binet_Padovan}, the equation \eqref{pleq:Pell_2} implies that
\begin{subequations}
\begin{align}\label{pleq:maj_n}
{\delta}^{n-3}\leq L_{n-2} \leq L_{n}-L_{n_1}=\mathcal{P}_m-\mathcal{P}_{m_1}<\alpha^{m-1},\\\label{pleq:min_n}
{\delta}^{n+1}\geq L_{n} > L_{n}-L_{n_1}=\mathcal{P}_m-\mathcal{P}_{m_1}\geq\mathcal{P}_{m-5}\geq \alpha^{m-7},
\end{align}
\end{subequations}
Hence
\begin{equation}\label{pleq:bound_n}
(m-7)\left(\dfrac{\log\alpha}{\log{\delta}}\right) -1 < n < (m-1)\left(\dfrac{\log\alpha}{\log{\delta}}\right) +3,
\end{equation}
where $\dfrac{\log\alpha}{\log{\delta}}=0.5843\ldots.$
If $n<$300, then $m\leq 190$. We ran a computer program for $2 \leq n_1 < n \leq 300$ and $1 \leq m_1 < m < 190$ and found only solutions from the list \eqref{pllist}. From now on we assume that $n \geq 300$.\\
Note that the inequality \eqref{pleq:bound_n} implies that $m <2n$. So, to solve the equation \eqref{pleq:Pell_2}, we need an upper bound on $n$.

\subsection{Upper bound on {\itshape n}}
Note that using the numerical inequalities\eqref{pleq:Numericle_alpha_beta_a_b} we have
\begin{equation}\label{plnumer_ineq}
	|\eta|^n +|\eta|^{n_1} + |b| |\beta|^m +|c| |\gamma|^m + |b| |\beta|^{m_1} +|c| |\gamma|^{m_1}<3.02.
\end{equation}
Using Binet's formulas in the Diophantine equation \eqref{pleq:Pell_2}, we get
$$
\begin{array}{rcl}
\abs {{\delta}^n-a\alpha^m} & = & \abs {-\eta^n+{\delta}^{n_1}+\eta^{n_1}+(b\beta^m+c\gamma^m) - (a\alpha^{m_1}+b\beta^{m_1}+c\gamma^{m_1})}\vspace{1mm}\\
  &\leq& {\delta}^{n_1}+ a \alpha^{m_1} +|\eta|^n + |\eta|^{n_1} + |b| |\beta|^m +|c| |\gamma|^m + |b| |\beta|^{m_1} +|c| |\gamma|^{m_1} \vspace{1mm}\\
  &<& {\delta}^{n_1}+ a \alpha^{m_1}+3.02 \vspace{1mm}\\
   &<& 4.76 \max\{{\delta}^{n_1},\alpha^{m_1}\}.
\end{array}
$$
By dividing by $a\alpha^m$ and using the relation \eqref{pleq:maj_n}, we obtain
$$
\begin{array}{lll}
\abs {a^{-1} {\delta}^n \alpha^{-m}-1}  &<&  \max\left\{ \dfrac{4.76}{a\alpha^m}{\delta}^{n_1},\dfrac{4.76}{a}\alpha^{m_1-m}\right\}\vspace{1mm}\\
  &<&  \max\left\{5.01 \dfrac{{\delta}^{n_1}}{\alpha \cdot {\delta}^{n-3}},5.01\alpha^{m_1-m}\right\}.
\end{array}
$$
Therefore, we get
\begin{equation}\label{pleq:Lambda}
\abs {a^{-1} {\delta}^n \alpha^{-m}-1}   < \max\{{\delta}^{n_1-n+7},\alpha^{m_1-m+7}\}.
\end{equation}
For the left member, we apply the theorem \ref{theoMatv} with the data
$$
s=3, \quad \gamma_1=a, \quad \gamma_2=\delta, \quad \gamma_3=\alpha, \quad b_1=-1, \quad b_2=n, \quad b_3=-m.
$$
Throughout our demonstrations, we work with $\KK:=\QQ(\sqrt{5},\alpha)$ with $D=6$. Since $\max\{1, n, m\}\leq 2n$ we take $B:=2n$. We have
$$
h(\gamma_2)=\dfrac{\log{\delta}}{2} \quad \text{and} \quad h(\gamma_3)=\dfrac{\log\alpha}{3}.
$$
Moreover, the minimal polynomial of $\gamma_1$ is $23x^3-23x^2+6x-1$ and has roots $a$, $b$, $c$. Since $|a|<1$ and $|b| = |c| < $1, then
$$
h(\gamma_1)=\dfrac{1}{3} \log23.
$$  
So we can take
$$
A_1=2\log23, \quad A_2=3\log{\delta}, \quad A_3= 2\log\alpha.
$$
Put
$$
\Lambda=a^{-1} {\delta}^n \alpha^{-m}-1.
$$
If $\Lambda=0$, then ${\delta}^n(\alpha^{-1})^m=a$, which is false, Since ${\delta}^n(\alpha^{ -1})^m\in\mathcal{O}_{\KK}$ while $a$ does not belong to $\mathcal{O}_{\KK}$, as can be seen immediately at from its minimal polynomial. Thus, $\Lambda\neq 0$. Then, by theorem \ref{theoMatv}, the left side of the equation \eqref{pleq:Lambda} is bounded by
$$
\log\abs\Lambda > -1.4\cdot 30^6 \cdot 3^{4.5} \cdot 6^2 (1+\log 6) (1+\log 2n) (2\log 23) (3\log{\delta}) (2\log\alpha).
$$
By comparing with \eqref{pleq:Lambda}, we get
$$
\min\{ (n-n_1-7) \log {\delta}, (m-m_1-7) \log\alpha \} < 7.33\times 10^{13} (1+\log 2n),
$$
which gives
$$
\min\{ (n-n_1) \log {\delta}, (m-m_1) \log\alpha \} < 7.33\times 10^{13} (1+\log 2n).
$$
Now two cases arise.\\
\textbf{Case 1.} $\min\{ (n-n_1) \log {\delta}, (m-m_1) \log\alpha \}= (n-n_1) \log {\delta}$.

\medskip

In this case, we rewrite \eqref{pleq:Pell_2} as
$$
\abs{\left({\delta}^{n-n_1}-1\right){\delta}^{n_1}-a\alpha^m} = \abs{-a\alpha^{m_1}+\eta^n- \eta^{n_1} +(b\beta^m+c\gamma^m) - (b\beta^{m_1}+c\gamma^{m_1})}
$$
using \eqref{plnumer_ineq} and dividing by $\alpha^m$, we get\begin{equation}\label{pleq:Lambda_1}
\abs {\left( \dfrac{{\delta}^{n-n_1}-1}{a} \right){\delta}^{n_1}\alpha^{-m} - 1} < 5.21 \alpha^{m_1-m}.
\end{equation}
We put
$$
\Lambda_1=\left( \dfrac{{\delta}^{n-n_1}-1}{a} \right){\delta}^{n_1}\alpha^{-m} - 1.
$$
We see that, $\Lambda_1\neq0$, for if $\Lambda_1=0$, then $ {\delta}^{n}-{\delta}^{n_1}=a\alpha^m$. This is impossible because $a\alpha^m\in \QQ(\alpha)$ and ${\delta}^{n}-{\delta}^{n_1}\in\QQ(\sqrt{5})\ backslash \QQ$. Indeed, if ${\delta}^{n}-{\delta}^{n_1}\in \QQ$, then when we take $\sigma \neq id$ to be the unique non-trivial $\QQ$- automorphism on $\QQ(\sqrt{5})$. Then we get
$$
{\delta}^n-{\delta}^{n_1}=\sigma({\delta}^n-{\delta}^{n_1})=\eta^n-\eta^{n_1}.
$$
However, the absolute value of the left member is at least ${\delta}^n-{\delta}^{n_1} \geq {\delta}^{n-2} \geq {\delta}^{298} >2$, while the absolute value on the right is at most $\abs{\eta^{n_1}-\eta^n} \leq |{\eta}|^{n_1}+|{\eta}|^ n <2.$ By this obvious contradiction, we conclude that $\Lambda_1\neq 0$.

We apply the lemma \ref{theoMatv} by taking $s=3$, and
$$
\gamma_1=\dfrac{{\delta}^{n-n_1}-1}{a}, \quad \gamma_2=\delta, \quad \gamma_3=\alpha, \quad b_1=1, \quad b_2=n_1, \quad b_3=-m.
$$
The minimal polynomial of ${\delta}^{n-n_1}+1$ divides
$$
x^2+(2-L_{n-n_1})x+((-1)^{n-n_1}+1-L_{n-n_1}),
$$
where $\{L_k\}_{k\geq 0}$ is the Lucas sequence defined by $L_0=2$, $L_1=1$, $L_{k+2}=L_{k+1}+L_k $ for all $k\geq0$, for which the Binet formula of its general term is
$$
L_k={\delta}^k+\eta^k\quad \text{for all }k\geq0.
$$

On the other hand, the minimal polynomial of $a$ is $23x^3-23x^2+6x-1$ and has roots $a$, $b$, $c$. Since $|b| = |c| < 1$ and $a<1$, then $h(a)=\dfrac{\log 23}{3}$.

Thus, we get
\begin{equation}\label{plmaj_h(gamma_1)}
\begin{array}{rcl}
h(\gamma_1) & \leq & h\left({\delta}^{n-n_1}+1 \right)+ h(a)\vspace{2mm}\\
  &\leq & h\left({\delta}^{n-n_1}\right)+h(1)+\log 2+ \dfrac{\log 23}{3}\vspace{2mm}\\
  &\leq&    (n-n_1)h\left({\delta}\right)+h(1)+\log 2+ \dfrac{\log 23}{3}\vspace{2mm}\\
  &<& \dfrac{1}{2} (n-n_1) \log \delta +1.74\vspace{2mm}\\
 &<&     3.66 \times 10^{13} (1+\log 2n).
\end{array}
\end{equation}
Thus, we can take $A_1:=2.2 \times 10^{14} (1+\log 2n)$. Also, as before, we can take $A_2:=3\log{\delta}$ and $A_3 :=2\log\alpha$. Finally, since $\max\{1,n_1,m\}\leq 2n $, we can take $B:=2n$.
We then obtain
$$
\log\abs{\Lambda_1} > -1.4\cdot 30^6 \cdot 3^{4.5} \cdot 6^2 (1+\log 6) (1+\log 2n) \times (2.2 \times 10^{14} (1+\log 2n)) (3\log{\delta}) (2\log\alpha).
$$
Thereby,
$$
\log\abs{\Lambda_1}> -2.57\cdot 10^{27} (1+\log 2n)^2.
$$
Comparing this with \eqref{pleq:Lambda_1}, we get that
$$
(m-m_1)\log \alpha<2.57\cdot 10^{28} (1+\log 2n)^2.
$$

\medskip
\textbf{Case 2.} $\min\{ (n-n_1) \log {\delta}, (m-m_1) \log\alpha \}= (m-m_1) \log\alpha $.

\medskip

In this case, we rewrite \eqref{pleq:Pell_2} as
$$
\abs{{\delta}^{n}-a\alpha^m+a\alpha^{m_1}} = \abs{\eta^n+ {\delta}^{n_1}-\eta^{n_1} +(b\beta^m+c\gamma^m) - (b\beta^{m_1}+c\gamma^{m_1})}
$$
then
\begin{equation}\label{pleq:Lambda_2}
\abs {\dfrac{{\delta}^n\alpha^{-m_1}}{ a (\alpha^{m-m_1}-1)} - 1} < \dfrac{4.03}{a(1-\alpha^{m_1-m})\alpha} \dfrac{{\delta}^{n_1}}{\alpha^{m-1}} < 
18 {\delta}^{n_1-n+3}.
\end{equation}
Let 
$$
\Lambda_2=(a(\alpha^{m-m_1}-1))^{-1} {\delta}^{n}\alpha^{-m_1} - 1.
$$
We see that, $\Lambda_2\neq0$, for if $\Lambda_2=0$ implies ${\delta}^{2n}=\alpha^{2m_1}a^2(\alpha^{m-m_1}-1 )^2.$, however, ${\delta}^{2n}\in\QQ(\sqrt{5})\backslash\QQ$, while $\alpha^{2m_1}a^2(\alpha^ {m-m_1}-1)^2\in\QQ(\alpha)$, which is not possible.
We apply again the lemma \ref{theoMatv}. In this application, we take $s=3$, and
$$
\gamma_1=a(\alpha^{m-m_1}-1), \quad \gamma_2=\delta, \quad \gamma_3=\alpha, \quad b_1=-1, \quad b_2=n, \quad b_3=-m_1.
$$
We have
$$
\begin{array}{lll}
h(\alpha^{m-m_1}-1)  &\leq& h(\alpha^{m-m_1})+h(-1)+\log 2=(m-m_1)h(\alpha)+\log 2\vspace{1mm}\\
 &=& \dfrac{(m-m_1)\log\alpha}{3} +\log 2 < 2.43 \times 10^{13} (1+\log 2n)+\log 2.
\end{array}
$$
Thus, on
$$
\begin{array}{lll}
h(\gamma_1) &<& 2.44 \times 10^{13} (1+\log 2n)+\log 2+\dfrac{\log 23}{3}+\log\sqrt{5}\vspace{1mm}\\
  &<& 2.44 \times 10^{13} (1+\log 2n).
\end{array}
$$
So we can take $A_1 :=1.47 \times 10^{14} (1+\log 2n)$. Also, as before, we can take $A_2:=3\log{\delta}$ and $A_3:=2\log\alpha$. Finally, since $\max\{1,n,m_1+1\}\leq 2n $, we can take $B:=2n$.
We then get this
$$
\log\abs{\Lambda_2} > -1.4\cdot 30^6 \cdot 3^{4.5} \cdot 6^2 (1+\log 6) (1+\log 2n) \times (1.47 \times 10^{14} (1+\log 2n)) (3\log{\delta}) (2\log\alpha).
$$
Hence,
$$
\log\abs{\Lambda_1}> -1.71\cdot 10^{27} (1+\log 2n)^2.
$$
Comparing this with \eqref{pleq:Lambda_2}, we get that
$$
(n-n_1)\log {\delta}<1.71\cdot 10^{27} (1+\log 2n)^2.
$$
Thus, in both cases 1 and 2, we have
\begin{subequations}
	\begin{align}\label{lmin-max}
	\min\{(n-n_1)\log{\delta},(m-m_1)\log\alpha\}<& 7.33\cdot 10^{13} (1+\log 2n)\\
	\max\{(n-n_1)\log{\delta},(m-m_1)\log\alpha\}<& 2.57\cdot 10^{27}(1+\log 2n)^2.\label{plmin-max}
	\end{align}
\end{subequations}

We finally rewrite the equation \eqref{pleq:Pell_2} as
$$
\abs{{\delta}^{n}-{\delta}^{n_1}-a\alpha^m+a\alpha^{m_1}} = \abs{\delta^n-\delta^{n_1} +(b\beta^m+c\gamma^m) - (b\beta^{m_1}+c\gamma^{m_1})} <3
$$
Dividing both sides by $a\alpha^{m_1}(\alpha^{m-m_1}-1)$, we get\begin{equation}\label{pleq:lambda_3}
\abs {\left( \dfrac{{\delta}^{n-n_1}-1}{a (\alpha^{m-m_1}-1)}\right) {\delta}^{n_1} \alpha^{-m_1}-1}< \dfrac{3}{a(1-\alpha^{m_1-m})\alpha} \dfrac{1}{\alpha^{m-1}}  < 13.5{\delta}^{3-n}. 
\end{equation}
To find a lower bound on the left side, we again use the lemma \ref{theoMatv} with $s=3$, and
$$
\gamma_1=\dfrac{{\delta}^{n-n_1}-1}{a (\alpha^{m-m_1}-1)}, \quad \gamma_2=\delta,\quad \gamma_3=\alpha,\quad b_1=1,\quad b_2=n_1, \quad b_3= -m_1.
$$
Using $h(x/y)=h(x)+h(y)$ for two nonzero algebraic numbers $x$ and $y$, we have
$$
\begin{array}{rcl}
h(\gamma_1) & \leq & h\left(\dfrac{{\delta}^{n-n_1}-1}{a}\right) +h(\alpha^{m-m_1}-1)\vspace{1mm}\\
  &<& \dfrac{1}{2}(n-n_1+4) \log {\delta} +\dfrac{\log 23}{3} + \dfrac{(m-m_1)\log\alpha}{3} +\log 2\vspace{1mm}\\
            & < & 2.14\cdot 10^{27} (1+\log 2n)^2,
\end{array}
$$
where in the chain of inequalities above, we used the argument of\eqref{plmaj_h(gamma_1)} as well as the \eqref{plmin-max} bound. So we can take $A_1:=1.54\cdot 10^{28} (1+\log 2n)^2$ and certainly $A_2:=3\log{\delta}$ and $A_3:=2\log\ alpha$. Using arguments similar to those in the proof that $\Lambda_1\neq 0$ we show that if we set
$$
\Lambda_3=\left( \dfrac{{\delta}^{n-n_1}-1}{a (\alpha^{m-m_1}-1)}\right) {\delta}^{n_1} \alpha^{-m_1}-1,
$$
then $\Lambda_3\neq 0$. The Lemma \ref{theoMatv} gives
$$
\log\abs{\Lambda_3} > -1.4\cdot 30^6 \cdot 3^{4.5} \cdot 6^2 (1+\log 6) (1+\log 2n) \times (1.54\cdot 10^{28} (1+\log 2n)^2) (3\log{\color{blue}{\delta}}) (2\log\alpha),
$$
which with \eqref{pleq:lambda_3} gives
$$
(n-3)<1.8\cdot 10^{41} (1+\log 2n)^3,
$$ 
leading to $n<2.45\cdot 10^{47}$.
\subsection{Reduction of the upper bound by {\itshape n}}
We must now reduce the bound above for $n$ and to do this we use the lemma \ref{lemDujella} several times and each time $M:=2.45\cdot 10^{47}$. To begin, let's go back to \eqref{pleq:Lambda} and set
$$
\Gamma:=n \log{\delta} -m \log \alpha-\log a.
$$
For technical reasons, we assume that $\min\{n-n_1,m-m_1\}\geq 20$. Let's go back to the inequalities for $\Lambda,~\Lambda_1,~\Lambda_2$.

Since we assume that $\min\{n-n_1,m-m_1\}\geq 20$, we get $|e^{\Gamma}-1|=|\Lambda|<\dfrac{1}{4} $. However, $|\Lambda|<\dfrac{1}{2}$ and since the inequality $|x|<2|e^x-1|$ holds for all $x\in\left(-\frac {1}{2}, \frac{1}{2}\right)$, we get
$$
\abs\Gamma < 2\max\{{\delta}^{n_1-n+7},\alpha^{m_1-m+7}\}\leq \max\{{\delta}^{n_1-n+9},\delta^{m_1-m+10}\}.
$$
Assume $\Gamma>0$. We then have the inequality
$$
\begin{array}{rcl}
0<n\left(\dfrac{\log{\delta}}{\log\alpha}\right)-m+\dfrac{\log(1/a)}{\log\alpha} &<& \max\left\{ \dfrac{{\delta}^{9}}{\log\alpha} {\delta}^{-(n-n_1)}, \dfrac{\alpha^{10}}{\log\alpha} \alpha^{-(m-m_1)} \right\} \vspace{2mm}\\
&<& \max\{ 270\cdot{\delta}^{-(n-n_1)},60\cdot \alpha^{-(m-m_1)} \}.
\end{array}
$$
We apply the lemma \ref{lemDujella} with
$$
\tau=\dfrac{\log{\delta}}{\log\alpha},\quad \mu=\dfrac{\log(1/a)}{\log\alpha},\quad (A,B)=(270,\delta)\; \text{ or }\; (60,\alpha).
$$
Let $\tau=[a_0,a_1,\ldots]=[1; 1, 2, 2, 6, 2, 1, 2, 1, 2, 1, 1, 11,\ldots]$ is the continued fraction of $\tau$. We consider the convergent $98$-th
$$
\dfrac{p}{q}=\dfrac{p_{98}}{q_{98}} = \dfrac{(78093067704223831799032754534503501859635391435517}{45634243076387457097046528084208490147594968308975}
$$
which satisfies $q=q_{98}>6M$. Moreover, this gives $\varepsilon>0.37$, and therefore either
$$
n-n_1\leq \dfrac{\log(270q/\varepsilon)}{\log{\delta}}<250,\;\text{ or }\;  m-m_1\leq \dfrac{\log(60q/\varepsilon)}{\log\alpha}< 420.
$$
In the case of $\Gamma<0$, we consider the following inequality:
$$
\begin{array}{rcl}
m\left(\dfrac{\log\alpha}{\log{\delta}}\right)-n+\dfrac{\log a}{\log{\delta}} &<& \max\left\{ \dfrac{{\delta}^{9}}{\log{\delta}} \alpha^{-(n-n_1)}, \dfrac{\alpha^{10}}{\log{\delta}} \alpha^{-(m-m_1)} \right\} \vspace{2mm}\\
&<& \max\{ 160\cdot{\delta}^{-(n-n_1)},37\cdot \alpha^{-(m-m_1)} \},
\end{array}
$$
instead and apply the lemma \ref{lemDujella} with
$$
\tau=\dfrac{\log\alpha}{\log{\delta}},\quad \mu=\dfrac{\log a}{\log{\delta}},\quad (A,B)=(160,\delta)\; \text{ or }\; (35,\alpha).
$$
Let $\tau=[a_0,a_1,\ldots]=[0,1,1,2,2,6,2,1,2,1,2,1,1,11,\ldots]$ be the fraction sequence of $\tau$ (note that the current $\tau$ is just the inverse of the previous $\tau$). Again, we consider the convergent $98$-th that satisfies $q=q_{98}>6M$. This again gives $\varepsilon>0.0867$, and so either
$$
n-n_1\leq \dfrac{\log(160q/\varepsilon)}{\log{\delta}}<246<250,\;\text{ or }\;  m-m_1\leq \dfrac{\log(35q/\varepsilon)}{\log\alpha}< 413<420.
$$
In conclusion, we have either $n-n_1\leq 250$ or $m-m_1\leq 420$ whenever $\Gamma\neq 0$.

Now, we must distinguish the cases $n-n_1\leq 250$ and $m-m_1\leq 420$. First suppose that $n-n_1\leq 250$. In this case, we consider the inequality \eqref{pleq:Lambda_1} and assume that $m-m_1\geq 20$. We ask
$$
\Gamma_1=n_1\log{\delta}-m\log\alpha+\log\left(\dfrac{{\delta}^{n-n_1}-1}{a}\right).
$$
Then the inequality \eqref{pleq:Lambda_1} implies that
$$
\abs{\Gamma_1} <10.4\alpha^{m_1-m}.
$$
If we further assume that $\Gamma_1>0$, then we get
$$
0<n_1\left(\dfrac{\log{\delta}}{\log\alpha}\right)-m+\dfrac{\log(({\delta}^{n-n_1}-1)/a)}{\log\alpha}<\dfrac{10.4}{(\log\alpha)}\alpha^{-(m-m_1)}<38\alpha^{-(m-m_1)}.
$$
We apply again the lemma \ref{lemDujella} with the same $\tau$ as in the case where $\Gamma>0$. We use the $100$-th ${p}/{q}={p_{98}}/{q_{98}}$ convergent to $\tau$ as before. But in this case we choose $(A,B):=(30,\alpha)$ and use
$$
\mu_k=\dfrac{\log(({\delta}^k-1)/a)}{\log\alpha},
$$
instead of $\mu$ for each possible value of $k:=n-n_1\in [1,2,\ldots 250].$ For the remaining values of $k$, we get $\varepsilon>0.00292$. Thus, according to the lemma \ref{lemDujella}, we obtain
$$
m-m_1< \dfrac{\log(38q/0.00292)}{\log\alpha}< 441.
$$ 
Thus, $n-n_1\leq 250$ implies $m-m_1\leq 441$.
In the case where $\Gamma_1 < 0$ we follow the ideas of the case where $\Gamma_1 > 0$. We use the same $\tau$ as in the case where $\Gamma < 0$ but instead of $\mu$ we take
$$
\mu_k=\dfrac{\log(a/({\delta}^k-1))}{\log{\delta}},
$$
for each possible value of $n- n_1 = k = 1, 2,\ldots , 250$. By using the lemma \ref{lemDujella} with this parameter we also obtain in this case that $n - n_1 \leq 250$ implies $m- m_1 \leq 435$.

In conclusion for $n - n_1 \leq 250$ we have $m- m_1 \leq 441$. Now let's go to the case where $m-m_1\leq 420$ and consider the inequality \eqref{pleq:Lambda_2}. we put
$$
\Gamma_2=n\log{\delta}-m_1\log\alpha+\log(1/(a(\alpha^{m-m_1}-1))),
$$
and we assume that $n-n_1 \geq 20$. We then have
$$
\abs{\Gamma_2} < \dfrac{36.1{\delta}^4}{\delta^{n-n_1}}.
$$
Assuming $\Gamma_2>0$, we get
$$
0<n\left(\dfrac{\log{\delta}}{\log\alpha}\right)-m_1+\dfrac{\log((1/(a(\alpha^{m-m_1}-1)))}{\log\alpha}<540\cdot{\delta}^{-(n-n_1)}.
$$
We apply the lemma \ref{lemDujella} again with the same $\tau$, $q$, $M$, $(A,B):=(540,\delta)$ and
$$
\mu_k=\dfrac{\log(1/(a(\alpha^{k}-1))}{\log \alpha}\quad \text{for }k=1, 2, \ldots 420.
$$
We obtain $\varepsilon>0.000354$, thus
$$
n-n_1<\dfrac{\log(540q/0.000354)}{\log{\delta}}<256.
$$
A similar conclusion is reached when $\Gamma_2<0$, indeed We get $\varepsilon>0.000508$, so
$$
n-n_1<\dfrac{\log(320q/0.000508)}{\log{\delta}}<256.
$$  

In conclusion, for $m-m_1\leq 420$ we have $n-n_1\leq 256$. So $m-m_1\leq 441$ and $n-n_1\leq 256$. Finally, we go to \eqref{pleq:lambda_3}. we put
$$
\Gamma_3=n_1\log{\delta}-m_1\log\alpha+\log\left(\dfrac{{\delta}^{n-n_1}-1}{a(\alpha^{m-m_1}-1)}\right).
$$
Since $n\geq 200$, the inequality \eqref{pleq:lambda_3} implies that
$$
\abs{\Gamma_3}<\dfrac{16.5}{{\delta}^{n-4}}=\dfrac{27{\delta}^4}{\delta^n}.
$$
Suppose $\Gamma_3>0$. Then 
$$
0<n_1\left(\dfrac{\log{\delta}}{\log\alpha}\right)-m_1+\dfrac{\log(({\delta}^k-1)/(a(\alpha^l-1)))}{\log\alpha}<240\cdot{\delta}^{-n},
$$
where $(k,l):=(n-n_1,m-m_1)$. We apply the lemma \ref{lemDujella} again with the same $\tau=\dfrac{\log{\delta}}{\log\alpha}$, $q_{98}$, $(A,B): =( 240,\delta)$ and
$$
\mu_{k,l}=\dfrac{\log(({\delta}^k-1)/(a(\alpha^l-1)))}{\log\alpha}\quad \text{for } 1\leq k\leq 256,\; 1\leq l\leq 441.
$$
We consider the $98$-th $\dfrac{p_{98}}{q_{98}}$ convergent. For all pairs $(k,l)$ we get that $\varepsilon>1.43\times 10^{-6}$. Thus, the lemma \ref{lemDujella} shows that
$$
n<\dfrac{\log (240 \times q_{98} \times 10^6/1.43)}{\log{\delta}}<271.
$$

The theorem \ref{plth:principal} is therefore proven.

\end{document}